\documentclass[reqno,9pt]{amsart}

\usepackage{comment}
\usepackage{indentfirst}
\usepackage[top=20mm, bottom=20mm, left=30mm, right=30mm]{geometry}
\usepackage{amscd}
\usepackage{amsmath,amsthm,amssymb}
\usepackage{color}
\usepackage[pdftex]{graphicx}
\usepackage{amsfonts}
\usepackage{cases}
\usepackage{bm}
\usepackage{pb-diagram}
\usepackage{braket}
\usepackage{multirow}
\usepackage{empheq}
\usepackage{tikz}

\usetikzlibrary{intersections,calc,arrows}

\numberwithin{equation}{subsection}
\theoremstyle{definition}
\newtheorem{defi}{Definition}[section]

\newtheorem{exa}[defi]{Example}

\newtheorem{theo}[defi]{Theorem}

\newcommand{\Ker}{\mathrm{Ker}}

\newcommand{\im}{\mathrm{Im}}

\newcommand{\initial}{\mathrm{in}}

\newcommand{\join}{\vee}
\newcommand{\meet}{\wedge}

\title{
Total of each of the 1-th and 2-th Betti numbers of the join-meet ideal of Special Crystal Lattice
}

\author{Yohei Oshida}

\address{Graduate School of Engineering and Science, Shibaura Institute of Technology, 307 Minumaku Fukasaku, Saitama-City 337-8570, Japan.}

\email{mf20022@shibaura-it.ac.jp}


\begin{document}

\begin{abstract}
\large
The join-meet ideal was introduced by Takayuki Hibi in 1987. It is binomial ideals that are defined by finite lattices. We study the join-meet ideal of non-distributive finite lattices that do not always satisfy modular. In particular, we work on the case of Crystal lattice which is one of them. In this paper, we give the computational result of a total of each of the 1-th and 2-th of grade Betti numbers of the join-meet ideal of a special Crystal lattice. The important point about this result is that it is characterized by the number of combinations of special generators. Moreover, we can consider from this result that a compatible monomial order suggests a connection to something good.
\end{abstract}

\maketitle{}

\setcounter{tocdepth}{1}
\tableofcontents

\section{Introduction}
Let $L$ be a finite lattice and $K[L]$ be the polynomial ring over a field $K$ whose variables are the elements of $L$. The ideal
\begin{eqnarray}
I_L=(\{ a b - (a \join b) (a \meet b) \mid a , b \in L \}) \subset K[L] \nonumber
\end{eqnarray}
is called the join-meet ideal of $L$. It was introduced in \cite{5} by Hibi. As shown by \cite{1} or \cite{5}, $L$ is distributive if and only if $I_{L}$ is a prime ideal. If $L$ is distributive lattices, $K[L]$ is called Hibi ring.

We study the join-meet ideal of non-distributive finite lattices that do not always satisfy modular. In particular, we work on the case of Crystal lattice, that is, the finite lattice $L_{k}(n_1,\cdots,n_ k)$ which is defined by
\begin{eqnarray}
L_{k}(n_1,\cdots,n_ k)=\{ s , x_{1,1} , \cdots , x_{1,n_1} , \cdots , x_{k,1} , \cdots , x_{k,n_k} , t \} \nonumber 
\end{eqnarray}
where $s < x_{i,1} < \cdots < x_{i,n_i} < t$ for $1 \leq i \leq k$. The definition of $L_{k}(n_1,\cdots,n_k)$ was introduced  by Yohei Oshida \cite{2} or \cite{4} in 2022. Note that this paper is the first to define its name as a crystal lattice.
We will not mention it at all here, but there is an important conjecture ''Crystal Conjecture'' with relation to the join-meet ideal of Crystal lattice; see \cite{2} or \cite{4} in detail.

Now, by using inverse lexicographic order induced by
\begin{eqnarray}
s < x_{1,1} < x_{1,2} < \cdots < x_{1,n_1} < x_{2,1} < x_{2,2} < \cdots < x_{2,n_2} <  \cdots < x_{k,1} < x_{k,2} < \cdots < x_{k,n_k} < t, \nonumber
\end{eqnarray}
we have
\begin{eqnarray}
\initial_{<}(I_{L_{k}(n_1,\cdots,n_k)})=( \bigcup_{i=1}^k \bigcup_{j=1,j \neq i}^k \{ x_{i,r_1} x_{j,r_2} \mid 1 \leq r_1 \leq n_i , 1 \leq r_2 \leq n_j \} ). \nonumber
\end{eqnarray}
Note that $<$ is a compatible monomial order; see [1, Example 6.16] in detail.

We recall grade Betti number by \cite{1} and \cite{3}. Let $R_{n_1,\cdots,n_k}=K[L_{k}(n_1,\cdots,n_k)]$. Let denote a minimal free resolution $S_{n_1,\cdots,n_k}$ of $R_{n_1,\cdots,n_k}$-module $R_{n_1,\cdots,n_k}/\initial_{<}(I_{L_{k}(n_1,\cdots,n_k)})$ by
\begin{eqnarray}
&S_{n_1,\cdots,n_k}:&
0
\leftarrow R_{n_1,\cdots,n_k}
\xleftarrow{\phi^{(n_1,\cdots,n_k)}_{1}} 
F^{(n_1,\cdots,n_k)}_1
\xleftarrow{\phi^{(n_1,\cdots,n_k)}_{2}} 
F^{(n_1,\cdots,n_k)}_2
\xleftarrow{\phi^{(n_1,\cdots,n_k)}_{3}}
\cdots
\nonumber \\
&&
\cdots
\xleftarrow{\phi^{(n_1,\cdots,n_k)}_{\ell_{n_1,\cdots,n_k}-1}}
F^{(n_1,\cdots,n_k)}_{\ell_{n_1,\cdots,n_k}-1}
\xleftarrow{\phi^{(n_1,\cdots,n_k)}_{\ell_{n_1,\cdots,n_k}}}
F^{(n_1,\cdots,n_k)}_{\ell_{n_1,\cdots,n_k}}
\leftarrow
0,
\nonumber
\end{eqnarray}
where $F^{(n_1,\cdots,n_k)}_i$ is a finitely generated grade free $R_{n_1,\cdots,n_k}$-module for $1 \leq i \leq \ell_{n_1,\cdots,n_k}$. Then, for module $F^{(n_1,\cdots,n_k)}_i$, we can write
\begin{eqnarray}
F^{(n_1,\cdots,n_k)}_i=\bigoplus_{j} R_{n_1,\cdots,n_k}(-X^{\bm{a}_{i,j}}), \nonumber
\end{eqnarray}
where $X=s x_{1,1} \cdots x_{1,n_1} x_{2,1} \cdots x_{2,n_2} \cdots x_{k,1} \cdots x_{k,n_k} t$ and $\bm{a}_{i,j}$ is the element of $\mathbb{Z}^{n_1+\cdots+n_k+2}$. Hence, since we denote $R_{n_1,\cdots,n_k}(-X^{\bm{a}_{i,j}})=R_{n_1,\cdots,n_k}(- |\bm{a}_{i,j}|)$ where $| \cdot |$ is $1$-norm, the minimal free resolution $S_{n_1,\cdots,n_k}$ can be rewritten as follows:
\begin{eqnarray}\label{minimal free resolution}
&S_{n_1,\cdots,n_k}:&
0
\leftarrow R_{n_1,\cdots,n_k}
\leftarrow
\bigoplus_{j} R_{n_1,\cdots,n_k}(- |\bm{a}_{1,j}|)^{B_{1,j}(n_1,\cdots,n_k)}
\leftarrow
\bigoplus_{j} R_{n_1,\cdots,n_k}(- |\bm{a}_{2,j}|)^{B_{2,j}(n_1,\cdots,n_k)}
\leftarrow
\nonumber \\
&&
\cdots 
\leftarrow
\bigoplus_{j} R_{n_1,\cdots,n_k}(- |\bm{a}_{\ell_{n_1,\cdots,n_k}-1,j}|)^{B_{\ell_{n_1,\cdots,n_k}-1,j}(n_1,\cdots,n_k)}
\nonumber \\
&&
\quad
\leftarrow
\bigoplus_{j} R_{n_1,\cdots,n_k}(- |\bm{a}_{\ell_{n_1,\cdots,n_k},j}|)^{B_{\ell_{n_1,\cdots,n_k},j}(n_1,\cdots,n_k)}
\leftarrow
0.
\end{eqnarray}
Then, $B_{i,j}(n_1,\cdots,n_k)$ is called grade Betti number. We are mainly interested in computing grade Betti number $B_{i,j}(n_1,\cdots,n_k)$ of minimal free resolution (\ref{minimal free resolution}). However, this computational is not easy. In general, for the ideal $I$ of $K[L]$ and a non-negative integer $i$ with $i \geq 2$, it is very difficult that we compute perfectively the total of each of the i-th Betti numbers of such $K[L]/I$ that $I$ has enough generators. Moreover, this computational becomes increasingly very complex when the number of generators increases regularly. In addition, the results may not be regular. Hence, it is very difficult that we compute perfectively the total of each of the i-th Betti numbers of $K[L_{k}(n_1,\cdots,n_k)]/I_{L_{k}(n_1,\cdots,n_k)}$. On the other hand, for $k=2$ and $n_2=1,2$, we can compute the total of each of the 1-th and 2-th Betti numbers by using [2, Theorem 1.1] or [4, Theorem 1.1]:
\begin{eqnarray}
\sum_{i=0}^{\infty} B_{1,i}(n_1,n_2), \quad
\sum_{i=0}^{\infty} B_{2,i}(n_1,n_2). \nonumber
\end{eqnarray}
We obtain the following result.

\begin{theo}\label{theo1} \it For $n_1 \geq 2$, we have
\begin{eqnarray}
\label{11} \sum_{i=0}^{\infty} B_{1,i}(n_1,1)
&=&
B_{1,2}(n_1,1)+B_{1,3}(n_1,1)
=
2n_1-1, \\
\label{12} \sum_{i=0}^{\infty} B_{2,i}(n_1,1)
&=&
B_{2,3}(n_1,1)+B_{2,4}(n_1,1)
=
n_1(n_1-1),
\end{eqnarray}
where
\begin{eqnarray}\label{combinational1}
B_{1,2}(n_1,1)=n_1, \quad 
B_{1,3}(n_1,1)=n_1-1, \quad
B_{2,3}(n_1,1)={}_{n_1} C_2, \quad
B_{2,4}(n_1,1)={}_{n_1-1} C_2 + n_1-1.
\end{eqnarray}
\end{theo}

\begin{theo}\label{theo2}\it For $n_1 \geq 2$, we have
\begin{eqnarray}
\label{21} \sum_{i=0}^{\infty} B_{1,i}(n_1,2)
&=&
B_{1,2}(n_1,2)+B_{1,3}(n_1,2)
=
3n_1, \\
\label{22} \sum_{i=0}^{\infty} B_{2,i}(n_1,2)
&=&
B_{2,3}(n_1,2)+B_{2,4}(n_1,2)
=
n_1(n_1+1)-1,
\end{eqnarray}
where
\begin{eqnarray}\label{combinational2}
B_{1,2}(n_1,2)=2n_1, \quad 
B_{1,3}(n_1,2)=n_1, \quad
B_{2,3}(n_1,2)=2 {}_{n_1} C_2 + n_1, \quad
B_{2,4}(n_1,2)={}_{n_1-1} C_2 + 3n_1-2.
\end{eqnarray}
\end{theo}

As we can see by reading the proof of Theorem \ref{theo1} and \ref{theo2}, (\ref{combinational1}) and (\ref{combinational2}) are Characterized by the number of combinations of special generators.

The structure of this paper is as follows. In section $2$, we introduce preparation to prove Theorem \ref{theo1} and \ref{theo2}. In particular, we mainly introduce binary operation to proof Theorem \ref{theo1} and \ref{theo2}. In section $3$, we proof Theorem \ref{theo1}. In section $4$, we prove Theorem \ref{theo2}. In section $5$, we introduce computational result by using Theorem \ref{theo1} and \ref{theo2}.\\

Below, unless otherwise noted, let $n_1=N$ and $x_{2,i}=y_i$ for $1 \leq i \leq n_2$. Moreover, let the vector $e_i$ be
\begin{eqnarray}
e_i=
\begin{cases}
1, & i\text{-th}, \\
0, & \text{otherwise}.
\end{cases}
\nonumber
\end{eqnarray}

\section{Preparation to prove Theorem \ref{theo1} and \ref{theo2}}

In this paper, we introduce binary operation $\otimes$ and $\boxplus$ to proof Theorem \ref{theo1} and \ref{theo2}.

At first, we introduce the definition of $\otimes$. In general, $\otimes$ is mainly a tensor product. However, note that it is not a tensor product in this paper.

\begin{defi}
Let denote the vector $u$ and $v$ by
\begin{eqnarray}
u=
\begin{bmatrix}
u_1 & u_2 & \cdots & u_n
\end{bmatrix}^T
, \quad 
v=
\begin{bmatrix}
v_1 & v_2 & \cdots & v_n
\end{bmatrix}^T
, \nonumber
\end{eqnarray}
where $u_1,\cdots,u_n,v_1,\cdots,v_n$ are the element of $R_{N,M}$.
Then, we define $u \otimes v$ by
\begin{eqnarray}
u \otimes v=
\begin{bmatrix}
v_1 & v_2 & \cdots & v_n
\\
u_1 & u_2 & \cdots & u_n
\end{bmatrix}^T
. \nonumber
\end{eqnarray}
\end{defi}

\begin{exa}
Let denote the vector $u$ and $v$ by
\begin{eqnarray}
u=
\begin{bmatrix}
1 & 2 & 3
\end{bmatrix}^T
, \quad 
v=
\begin{bmatrix}
3 & 2 & 1
\end{bmatrix}^T
. \nonumber
\end{eqnarray}
Then, $u \otimes v$ is computed as follows:
\begin{eqnarray}
u \otimes v=
\begin{bmatrix}
1 & 3 \\
2 & 2 \\
3 & 1
\end{bmatrix}
. \nonumber
\end{eqnarray}
\end{exa}

Next, we introduce the definition of $\boxplus$.

\begin{defi}Let denote the vector $u$ and $v$ by
\begin{eqnarray}
u=
\begin{bmatrix}
u_1 & u_2 & \cdots & u_n
\end{bmatrix}^T
, \quad 
v=
\begin{bmatrix}
v_1 & v_2 & \cdots & v_m
\end{bmatrix}^T
, \nonumber
\end{eqnarray}
where $u_1,\cdots,u_n,v_1,\cdots,v_m$ are the element of $R_{N,M}$.
Then, we define $u \boxplus v$ by
\begin{eqnarray}
u \boxplus v=
\begin{bmatrix}
u_1 & u_2 & \cdots & u_n
&
v_1 & v_2 & \cdots & v_n
\end{bmatrix}^T
. \nonumber
\end{eqnarray}
\end{defi}

\begin{exa}
Let denote the vector $u$ and $v$ by
\begin{eqnarray}
u=
\begin{bmatrix}
1 & 2 & 3
\end{bmatrix}^T
, \quad 
v=
\begin{bmatrix}
2 & 11 & 4 & 21
\end{bmatrix}^T
. \nonumber
\end{eqnarray}
Then, $u \boxplus v$ is computed as follows:
\begin{eqnarray}
u \boxplus v=
\begin{bmatrix}
1 & 2 & 3 & 2 & 11 & 4 & 21
\end{bmatrix}^T
. \nonumber
\end{eqnarray}
\end{exa}

\section{The proof of Theorem \ref{theo1}}

At first, for $N=2$, we show that equation (\ref{11}) and (\ref{12}) hold.

\subsection{Construction of $F^{(2,1)}_1$ and computational of total Betti number of $F^{(2,1)}_1$}

By [2, Theorem 1.1] or [4, Theorem 1.1], we have
\begin{eqnarray}
\initial_{<}(I_{L_2(2,1)})=\{ x_i y_1 \mid 1 \leq i \leq 2 \} \cup \{ x_2 s t \}. \nonumber
\end{eqnarray}
Thus, since we have $\im(\phi^{(2,1)}_{1})=\initial_{<}(I_{L_2(2,1)})$, we obtain
\begin{eqnarray}
F^{(2,1)}_1
=\Biggr ( \bigoplus_{i=1}^2 R_{2,1}(x_i y) \Biggr ) \oplus \Biggr ( \bigoplus_{i=2}^2 R_{2,1}(x_i s t) \Biggr ). \nonumber
\end{eqnarray}
Then, free module $F^{(2,1)}_1$ can be written as follows:
\begin{eqnarray}
F^{(2,1)}_1
=
\Biggr ( \bigoplus_{i=1}^2 R_{2,1}(-2) \Biggr )
\oplus
\Biggr ( \bigoplus_{i=2}^2 R_{2,1}(-3) \Biggr )
=
\bigoplus_{i=1} R_{2,1}(-2)^2
\oplus
R_{2,1}(-3).
\nonumber
\end{eqnarray}
Hence, it follows from the above result that we have
\begin{eqnarray}\label{0000}
\sum_{i=0}^{\infty} B_{1,i}(2,1)=B_{1,2}(2,1)+B_{1,3}(2,2)=2+1=3.
\end{eqnarray}
Therefore, since the equation (\ref{0000}) equals (\ref{11}) with $N=2$, we showed that equation (\ref{11}) holds for $N=2$.\\

\subsection{Construction of $F^{(2,1)}_{2}$ and computational of total Betti number of $F^{(2,1)}_{2}$}

Since we have $\im(\phi^{(2,1)}_{2})=\Ker(\phi^{(2,1)}_{1})$, we have
\begin{eqnarray}
\im(\phi^{(2,1)}_{2})=\Biggr \{ u=\sum_{i=1}^2 a_i x_i y_1 e_i + b x_2 s t e_3 \mid a_i,b \in R_{2,1},\phi^{(2,1)}_{1}(u)=0 \Biggr \}. \nonumber
\end{eqnarray}
Then, the necessary and sufficient conditions to calculate $u$ satisfying
\begin{eqnarray}\label{001}
\phi^{(2,1)}_{1}(u)=\sum_{i=1}^2 a_i x_i y_1 + b x_2 s t=0
\end{eqnarray}
are to calculate $u$ satisfying the following each equation:
\begin{eqnarray}\label{001.1}
a_1 x_1 y_1 + a_2 x_2 y_1&=&0, \\
\label{001.2}
a_1 x_1 y_1 + b x_2 s t&=&0, \\
\label{001.3}
a_2 x_2 y_1 + b x_2 s t&=&0.
\end{eqnarray}
By solving (\ref{001.1}), we obtain
\begin{eqnarray}
u=\alpha ( e_1 - e_2 ) x_1 x_2 y_1 + \beta (e_1 - e_3) x_1 x_2 y_1 s t + \gamma ( e_2 - e_3 ) x_2 y_1 s t,  \nonumber
\end{eqnarray}
where $\alpha$, $\beta$, $\gamma$ is the element of $R_{2,1}$. We can rearrange the above equation by seeting $\beta=0$:
\begin{eqnarray}\label{001.4}
u
&=&
\alpha ( e_1 - e_2 ) x_1 x_2 y_1
+
\biggr ( \beta (e_1 - e_3) x_1 + \gamma ( e_2 - e_3 ) \biggr ) x_2 y_1 s t
\nonumber \\
&=&
\alpha ( e_1 - e_2 ) x_1 x_2 y_1 + \gamma ( e_2 - e_3 ) x_2 y_1 s t
\nonumber \\
&=&
\begin{bmatrix}
1 & 0 \\ 
-1 & 1 \\
0 & -1
\end{bmatrix}
\begin{bmatrix}
\alpha x_1 x_2 y_1 \\ \gamma x_2 y_1 s t 
\end{bmatrix}.
\end{eqnarray}
Thus, since we have (\ref{001.4}) and $\im(\phi^{(2,1)}_{3})=\Ker(\phi^{(2,1)}_{2})=0$, we have $F^{(2,1)}_{2}=R_{2,1}(x_1 x_2 y_1) \oplus R_{2,1}(x_2 y_1 s t)$. Then, since we have $F^{(2,1)}_{2}=R_{2,1}(-3) \oplus R_{2,1}(-4)$, we have
\begin{eqnarray}\label{001.5}
\sum_{i=0}^{\infty} B_{2,i}(2,1)=B_{2,2}(2,1)+B_{2,3}(2,1)=1+1=2.
\end{eqnarray}
Hence, since the equation (\ref{001.5}) equals (\ref{12}) with $N=2$, we showed that equation (\ref{12}) holds for $N=2$.\\

Next, for $N \geq 3$, we show that equation (\ref{11}) holds.

\subsection{Construction of $F^{(N,1)}_{1}$ and computational of total Betti number of $F^{(N,1)}_{1}$}

By [2, Theorem 1.1] or [4, Theorem 1.1], we have
\begin{eqnarray}
\initial_{<}(I_{L_2(N,1)})=( \{ x_i y_1 \mid 1 \leq i \leq N \} \cup \{ x_i s t \mid 2 \leq i \leq N \} ). \nonumber
\end{eqnarray}
Thus, it follows from $\im(\phi^{(N,1)}_{1})=\initial_{<}(I_{L_2(N,1)})$ that we obtain
\begin{eqnarray}
F^{(N,1)}_{1}
=\Biggr ( \bigoplus_{i=1}^N R_{N,1}(x_i y) \Biggr ) \oplus \Biggr ( \bigoplus_{i=2}^N R_{N,1}(x_i s t) \Biggr ). \nonumber
\end{eqnarray}
Then, free module $F^{(N,1)}_{1}$ can be written as follows:
\begin{eqnarray}
F^{(N,1)}_{1}
=
\Biggr ( \bigoplus_{i=1}^N R_{N,1}(-2) \Biggr )
\oplus
\Biggr ( \bigoplus_{i=2}^N R_{N,1}(-3) \Biggr )
=
\bigoplus_{i=1} R_{N,1}(-2)^{N}
\oplus
R_{N,1}(-3)^{N-1}.
\nonumber
\end{eqnarray}
Hence, it follows from the above result that we have
\begin{eqnarray}
\sum_{i=0}^{\infty} B_{1,i}(N,1)=B_{1,2}(N,1)+B_{1,3}(N,1)=N+(N-1)=2N-1. \nonumber
\end{eqnarray}
Therefore, we showed that equation (\ref{11}) holds.\\

Finally, we show that equation (\ref{12}) holds.

\subsection{The computational of $\im(\phi^{(N,1)}_{2})$} \label{sub1}

By $\im(\phi^{(N,1)}_{2})=\Ker(\phi^{(N,1)}_{1})$, we have
\begin{eqnarray}
\im(\phi^{(N,1)}_{2})=\Biggr \{ u=\sum_{i=1}^N a_i x_i y_1 e_i + \sum_{i=2}^N b_i x_i s t e_{N+i-1} ; a_i, b_i \in R_{N,1}, \phi^{(N,1)}_{1}(u)=0 \Biggr \}.
\nonumber
\end{eqnarray}
Then, the necessary and sufficient conditions to calculate $u$ satisfying
\begin{eqnarray}\label{31}
\phi^{(N,1)}_{1}(u)=\sum_{i=1}^N a_i x_i y_1 + \sum_{i=2}^N b_i x_i s t=0
\end{eqnarray}
are to calculate $u$ satisfying $\phi^{(N,1)}_{1}(u)=0$ with the following cases:
\begin{eqnarray}
&&\text{$1$. The case of $b_i=0$} \nonumber \\
&&\text{$2$. The case of $a_i=0$} \nonumber \\
&&\text{$3$. The case of $a_i \neq 0, b_i \neq 0$} \nonumber
\end{eqnarray}
Below, we compute $u$ satisfying $(\ref{31})$ with each of the above cases.

\subsubsection{\textup{The case of $b_i=0$}} In this case, since we have $(\ref{31})$, we have
\begin{eqnarray}\label{31.1}
\sum_{i=1}^{N} a_i x_i y_1= 0. 
\end{eqnarray}
The necessary and sufficient conditions to satisfy $(\ref{31.1})$ are to satisfy
\begin{eqnarray}\label{31.12}
a_i x_i y_1 + a_j x_j y_1 = 0
\end{eqnarray}
for $1 \leq i < j \leq N$. By solving (\ref{31.12}), we obtain
\begin{eqnarray}
(a_i,a_j)=\alpha_{i+j}(x_j,- x_i), \nonumber
\end{eqnarray}
where $\alpha_{i+j}$ is the element of $R_{N,1}$. Hence, the vector $u$ that satisfies $(\ref{31.1})$ is
\begin{eqnarray}\label{31.13}
u_1=\sum_{i=1}^{N} \sum_{j=1,j > i}^{N} \alpha_{i+j} ( x_i x_j y_1 e_i - x_i x_j y_1 e_j).
\end{eqnarray}

(The case of $a_i=0$) In this case, since we have $(\ref{31})$, we have
\begin{eqnarray}\label{31.2}
\sum_{i=2}^N b_i x_i s t=0.
\end{eqnarray}
The necessary and sufficient conditions to satisfy $(\ref{31.2})$ are to satisfy
\begin{eqnarray}\label{31.22}
b_i x_i s t + b_j x_j s t = 0 
\end{eqnarray}
for $2 \leq i < j \leq N$. By solving (\ref{31.22}), we obtain
\begin{eqnarray}
(b_i,b_j)=\beta_{i+j}(x_j,-x_i) \nonumber
\end{eqnarray}
where $\beta_{i+j}$ is the element of $R_{N,1}$. Hence, the vector $u$ that satisfies $(\ref{31.2})$ is
\begin{eqnarray}\label{31.23}
u_2=\sum_{i=2}^{N} \sum_{j=3,i < j}^N \beta_{i+j} ( x_i x_j s t e_{i+N-1} - x_i x_j s t e_{j+N-1} ).
\end{eqnarray}

(The case of $a_i \neq 0, b_i \neq 0$) In this case, since we have $(\ref{31.13})$ and $(\ref{31.23})$, we only have to compute $(a_i,b_j)$ satisfying
\begin{eqnarray}\label{31.3}
a_i x_i y_1 + b_j x_j s t = 0
\end{eqnarray}
for $1 \leq i \leq N$ and $2 \leq j \leq N$.
The computational result is as follows:
\begin{eqnarray}
(a_i,b_j)=
\begin{cases}
\gamma_{2i} (s t , - y_1), & i=j, \\
\gamma_{i+j} (x_j s t , - x_i y_1), & i \neq j,
\end{cases}
\nonumber
\end{eqnarray}
where $\gamma_{i+j}$ is the polynomial of $R_{N,1}$. Hence, the vector $u$ that satisfies $(\ref{31.3})$ is 
\begin{eqnarray}\label{31.31}
u_3=
\sum_{i=2}^N \gamma_{2i} ( x_i y_1 s t e_i - x_i y_1 s t e_{N+i-1} )
+
\sum_{i=1}^N \sum_{j=2,j \neq i}^N \gamma_{i+j} ( x_i x_j y_1 s t e_i - x_i x_j y_1 s t e_{N+j-1} ).
\end{eqnarray}
By setting $\gamma_{i+j}=0$ for $i \neq j$, We can rearrange (\ref{31.31}) as follows:
\begin{eqnarray}
&&
\sum_{i=2}^N \gamma_{2i} ( x_i y_1 s t e_i - x_i y_1 s t e_{N+i-1} )
+
\sum_{i=1}^N \sum_{j=2,j \neq i}^N \gamma_{i+j} ( x_i x_j y_1 s t e_i - x_i x_j y_1 s t e_{N+j-1} )
\nonumber \\
&=&
\sum_{i=2}^N \gamma_{2i} ( x_i y_1 s t e_i - x_i y_1 s t e_{N+i-1} )
+
\sum_{j=2}^N \gamma_{1+j} ( x_1 x_j y_1 s t e_1 - x_1 x_j y_1 s t e_{N+j-1} )
\nonumber \\
&&
+
\sum_{i=2}^N \sum_{j=2,j \neq i}^N \gamma_{i+j} ( x_i x_j y_1 s t e_i - x_i x_j y_1 s t e_{N+j-1} )
\nonumber \\
&=&
\sum_{i=2}^N \gamma_{2i} x_i y_1 s t e_i
-
\sum_{i=2}^N ( \gamma_{2i} + \gamma_{1+i} x_1 ) x_i y_1 s t e_{N+i-1} 
+
\sum_{i=2}^N \gamma_{1+i} x_1 x_i y_1 s t e_1 
\nonumber \\
&&
+
\sum_{i=2}^N \sum_{j=2,j \neq i}^N \gamma_{i+j} ( x_i x_j y_1 s t e_i - x_i x_j y_1 s t e_{N+j-1} )
\nonumber \\
&=&
\sum_{i=2}^N \gamma_{2i} ( x_i y_1 s t e_i - x_i y_1 s t e_{N+i-1} )
+
\sum_{i=2}^N \sum_{j=2,j \neq i}^N \gamma_{i+j} ( x_i x_j y_1 s t e_i - x_i x_j y_1 s t e_{N+j-1} )
\nonumber \\
&=&
\sum_{i=2}^N
\Biggr (
\gamma_{2i}
+ 
\sum_{j=2,j \neq i}^N \gamma_{i+j} x_j 
\Biggr ) x_i y_1 s t e_i
-
\sum_{i=2}^N 
\Biggr (
\gamma_{2i} x_i y_1 s t e_{N+i-1}
+
\sum_{j=2,j \neq i}^N \gamma_{i+j} x_i x_j y_1 s t e_{N+j-1}
\Biggr )
\nonumber \\
&=&
\sum_{i=2}^N \gamma_{2i} ( x_i y_1 s t e_i - x_i y_1 s t e_{N+i-1} ).
\nonumber
\end{eqnarray}
Therefore, we have
\begin{eqnarray}\label{31.32}
u_3=
\sum_{i=2}^N \gamma_{2i} ( x_i y_1 s t e_i - x_i y_1 s t e_{N+i-1} )
.
\end{eqnarray} 

\subsection{Matrix representation of $u$}
By (\ref{31.13}), (\ref{31.23}) and (\ref{31.32}), we have $u=\sum_{i=1}^3 u_i$. We rearrange $u$ to emphasize that there is no negative component in the 1-th component:
\begin{eqnarray}\label{31.4}
u
&=&
\sum_{j=2}^{N} \alpha_{1+j} ( x_1 x_j y_1 e_1 - x_1 x_j y_1 e_j)
+
\sum_{i=2}^{N} \sum_{j=1,j > i}^{N} \alpha_{i+j} ( x_i x_j y_1 e_i - x_i x_j y_1 e_j)
\nonumber \\
&&
+
\sum_{i=2}^{N} \sum_{j=3,i < j}^N \beta_{i+j} ( x_i x_j s t e_{i+N-1} - x_i x_j s t e_{j+N-1} )
+
\sum_{i=2}^N \gamma_{2i} ( x_i y_1 s t e_i - x_i y_1 s t e_{N+i-1} ).
\end{eqnarray}
Moreover, by using new vector $d_{i,j}=e_i - e_j$ for $1 \leq i , j \leq 2N-1$, we can rearrange (\ref{31.4}) as follows:
 \begin{eqnarray}\label{31.5}
u
&=&
\sum_{j=2}^{N} \alpha_{1+j} ( e_1 - e_j) x_1 x_j y_1
+
\sum_{i=2}^{N} \sum_{j=3,j > i}^{N} \alpha_{i+j} ( e_i - e_j) x_i x_j y_1
\nonumber \\
&&
+
\sum_{i=2}^{N} \sum_{j=3,j > i}^N \beta_{i+j} ( e_{i+N-1} - e_{j+N-1} ) x_i x_j s t
+
\sum_{i=2}^N \gamma_{2i} ( e_i - e_{N+i-1} ) x_i y_1 s t
\nonumber \\
&=&
\sum_{j=2}^{N} \alpha_{1+j} d_{1,j}  x_1 x_j y_1
+
\sum_{i=2}^{N} \sum_{j=3,j > i}^{N} \alpha_{i+j} d_{i,j} x_i x_j y_1
\nonumber \\
&&
+
\sum_{i=2}^{N} \sum_{j=3,j > i}^N \beta_{i+j} d_{i+N-1,j+N-1} x_i x_j s t
+
\sum_{i=2}^N \gamma_{2i} d_{i,N+i-1} x_i y_1 s t.
\end{eqnarray}
We introduce the symbols necessary to express (\ref{31.5}) by using a matrix. Let denote $f_{N,1}$ and $g_{N,1}$ by
\begin{eqnarray}
f_{N,1}
&=&
\Biggr ( \bigotimes_{i=2}^N d_{1,i} \Biggr ) 
\otimes
\Biggr ( \bigotimes_{i=2}^N \bigotimes_{j=3,j > i}^N d_{i,j} \Biggr )
\otimes \Biggr ( \bigotimes_{i=2}^N \bigotimes_{j=3,j > i}^N d_{i+N-1,j+N-1} \Biggr )
\otimes
\Biggr ( \bigotimes_{i=2}^N d_{i,N+i-1} \Biggr ),
\nonumber \\
g_{N,1}
&=&
v_{1,1} \boxplus v_{1,2} \boxplus \cdots \boxplus v_{1,N-1} \boxplus v_{2,1} \boxplus v_{2,2} \boxplus \cdots \boxplus v_{2,N-1} \boxplus v_{3,N}
\nonumber
\end{eqnarray}
where
\begin{eqnarray}
v_{1,i}
&=&
\begin{pmatrix}
x_i x_{i+1} y_1 & x_i x_{i+2} y_1 & \cdots & x_i x_N y_1 \\
\end{pmatrix}
, \quad \text{for $1 \leq i \leq N-1$},
\nonumber \\
v_{2,i}
&=&
\begin{pmatrix}
x_i x_{i+1} s t & x_i x_{i+2} s t & \cdots & x_i x_N s t \\
\end{pmatrix}
, \quad \text{for $2 \leq i \leq N-1$},
\nonumber \\
v_{3,N}
&=&
\begin{pmatrix}
x_2 y_1 s t & x_3 y_1 s t & \cdots & x_N y_1 s t \\
\end{pmatrix}
.
\nonumber
\end{eqnarray}
By using $f_{N,1}$ and $g_{N,1}$, it follows from (\ref{31.5}) that we have $u=f_{N,1}{g_{N,1}}^T$. Hence, we have a free resolution
\begin{eqnarray}\label{resolution_1}
0 \leftarrow R_{N,1} \xleftarrow{
\overbrace{
\begin{bmatrix}
1 & \cdots & 1
\end{bmatrix}
}^{2N-1 \text{pieces}}
} F^{(N,1)}_{1} \xleftarrow{f_{N,1}} F'_{2} \leftarrow \cdots
\end{eqnarray}
where
\begin{eqnarray}
F'_{2}=
\Biggr ( \bigoplus_{i=1}^{N-1} \bigoplus_{j=2,i<j}^N R_{N,1}(x_i x_j y_1) \Biggr )
\oplus 
\Biggr ( \bigoplus_{i=2}^{N-1} \bigoplus_{j=3,i<j}^N R_{N,1}(x_i x_j s t) \Biggr )
\oplus 
\Biggr ( \bigoplus_{i=2}^{N} R_{N,1}(x_i y_1 s t) \Biggr )
.
\end{eqnarray}

\subsection{Construction of $F^{(N,1)}_{2}$}
Then, we show that $F^{(N,1)}_{2}=F'_{2}$. Suppose that $F^{(N,1)}_{2} \neq F'_{2}$. Since we have $F^{(N,1)}_{2} \subsetneq F'_{2}$, for instance we can suppose that
\begin{eqnarray}
F^{(N,1)}_{2}
=
\Biggr ( \bigoplus_{i=1}^{N-1} \bigoplus_{j=2,i<j}^N R_{N,1}(x_i x_j y_1) \Biggr ) \oplus 
\Biggr ( \bigoplus_{i=2}^{N} R_{N,1}(x_i y_1 s t) \Biggr ).
\nonumber
\end{eqnarray}
However, such assumption means that we have $F^{(N,1)}_{2}=0$. In fact, if we have $F^{(N,1)}_{2} \neq 0$, then we have $\im(\phi^{(N,1)}_{2}) \subsetneq \text{span}_{R_{N,1}} \{ \omega \mid \text{$\omega$ is each term of $u=\sum_{i=1}^3 u_i$.} \}=\Ker(\phi^{(N,1)}_{1})$. It contradicts that we have $\im(\phi^{(N,1)}_{2})=\Ker(\phi^{(N,1)}_{1})$. Since we have $F^{(N,1)}_{2}=0$, it contradicts existence of the free resolution (\ref{resolution_1}). Hence, we showed that $F^{(N,1)}_{2}=F'_{2}$.

\subsection{Computational of total Betti number of $F^{(N,1)}_{2}$}
By $F^{(N,1)}_{2}=F'_{2}$, we have
\begin{eqnarray}
F^{(N,1)}_{2}
=
\Biggr ( \bigoplus_{i=1}^{N-1} \bigoplus_{j=2,i<j}^N R_{N,1}(x_i x_j y_1) \Biggr )
\oplus 
\Biggr ( \bigoplus_{i=2}^{N-1} \bigoplus_{j=3,i<j}^N R_{N,1}(x_i x_j s t) \Biggr )
\oplus 
\Biggr ( \bigoplus_{i=2}^{N} R_{N,1}(x_i y_1 s t) \Biggr )
. \nonumber
\end{eqnarray}
The above free module $F^{(N,1)}_{2}$ can be written as follows:
\begin{eqnarray}
F^{(N,1)}_{2}
&=&
R_{N,1}(-3)^{{}_N C_2}
\oplus
R_{N,1}(-4)^{{}_{N-1} C_2}
\oplus
R_{N,1}(-4)^{N-1}
\nonumber \\
&=&
R_{N,1}(-3)^{{}_N C_2}
\oplus
R_{N,1}(-4)^{{}_{N-1} C_2 + N-1}
.
\nonumber
\end{eqnarray}
Hence, it follows from the above result that we have
\begin{eqnarray}
\sum_{i=0}^{\infty} B_{2,i}(N,1)
&=&
B_{2,3}(N,1)
+
B_{2,4}(N,1)
\nonumber \\
&=&
{}_N C_2 + ({}_{N-1} C_2 + N-1)
\nonumber \\
&=&
\frac{N(N-1)}{2}+\frac{(N-1)(N-2)}{2}+\frac{2(N-1)}{2}
\nonumber \\
&=&
\frac{(N^2-N)+(N^2-3N+2)+(2N-2)}{2}
\nonumber \\
&=&\frac{2 N^2 - 2 N}{2}
\nonumber \\
&=&N(N-1).
\nonumber
\end{eqnarray}
Therefore, we showed that equation (\ref{12}) holds.

\section{The proof of Theorem \ref{theo2}}

At first, we show that equation (\ref{21}) holds.

\subsection{Calculation of $F^{(N,2)}_{1}$ and computational of total Betti number of $F^{(N,2)}_{1}$}

By [2, Theorem 1.1] or [4, Theorem 1.1], we have
\begin{eqnarray}
\initial_{<}(I_{L_2(N,2)})=( \{ x_i y_j \mid 1 \leq i \leq N , 1 \leq j \leq 2 \} \cup \{ x_i s t \mid 2 \leq i \leq N \} \cup \{ y_2 s t \} ). \nonumber
\end{eqnarray}
Thus, it follows from $\im(\phi^{(N,2)}_{1})=\initial_{<}(I_{L_2(N,2)})$ that we have
\begin{eqnarray}
F^{(N,2)}_1=\Biggr ( \bigoplus_{i=1}^N \bigoplus_{j=1}^2 R_{N,2}(x_i y_j) \Biggr ) \oplus \Biggr ( \bigoplus_{i=2}^N R_{N,2}(x_i s t) \Biggr ) \oplus R_{N,2}(y_2 s t). \nonumber
\end{eqnarray}
Then, free module $F^{(N,2)}_1$ can be written as follows:
\begin{eqnarray}
F^{(N,2)}_1
&=&
\Biggr ( \bigoplus_{i=1}^N \bigoplus_{j=1}^2 R_{N,2}(-2) \Biggr )
\oplus
\Biggr ( \bigoplus_{i=2}^N R_{N,2}(-3) \Biggr )
\oplus R_{N,2}(-3)
\nonumber \\
&=&
\bigoplus_{j=1}^2 R_{N,2}(-2)^{2N}
\oplus
R_{N,2}(-3)^N.
\nonumber
\end{eqnarray}
Hence, it follows from the above result that we have
\begin{eqnarray}
\label{theo21} \sum_{i=0}^{\infty} B_{1,i}(N,2)
&=&
B_{1,2}(N,2)+B_{1,3}(N,2)
\nonumber \\
&=&
2N+N
\nonumber \\
&=&3N.
\nonumber
\end{eqnarray}
Therefore, we showed that equation (\ref{21}) holds.\\

Next, we show that equation (\ref{22}) holds.

\subsection{The computational of $\im(\phi^{(N,2)}_{2})$} \label{sub1}

By $\im(\phi^{(N,2)}_{2})=\Ker(\phi^{(N,2)}_{1})$, we have
\begin{eqnarray}
&&
\im(\phi^{(N,2)}_{2})
\nonumber \\
&=&\Biggr \{ u=\sum_{i=1}^{N} a_{1,i} x_i y_1 e_i + \sum_{i=1}^{N} a_{2,i} x_i y_2 e_{N+i} + \sum_{i=2}^{N} b_i x_i s t e_{2N+i-1} + c y_2 s t e_{3N} ; a_{1,i},a_{2,i},b_i,c \in R_{N,2}, \phi^{(N,2)}_{1}(u)=0 \Biggr \}
\nonumber
\end{eqnarray}
Then, the necessary and sufficient conditions to calculate $u$ satisfying
\begin{eqnarray}\label{41}
\phi^{(N,2)}_{1}(u)=\sum_{i=1}^{N} a_{1,i} x_i y_1 + \sum_{i=1}^{N} a_{2,i} x_i y_2 + \sum_{i=2}^{N} b_i x_i s t  + c y_2 s t = 0 
\end{eqnarray}
are to calculate $u$ satisfying $\phi^{(N,2)}_{1}(u)=0$ with the following cases:
\begin{eqnarray}
&&\text{$1$. The cases of $a_{2,i},b_i,c=0$} \nonumber \\
&&\text{$2$. The cases of $a_{1,i},b_i,c=0$} \nonumber \\
&&\text{$3$. The cases of $a_{1,i},a_{2,i},c=0$} \nonumber \\
&&\text{$4$. The cases of $b_i,c=0$} \nonumber \\
&&\text{$5$. The cases of $a_{2,i},c=0$} \nonumber \\
&&\text{$6$. The cases of $a_{2,i},b_i=0,c \neq 0$} \nonumber \\
&&\text{$7$. The cases of $a_{1,i},c=0$} \nonumber \\
&&\text{$8$. The cases of $a_{1,i},b_i=0,c \neq 0$} \nonumber \\
&&\text{$9$. The cases of $a_{1,i},a_{2,i}=0,c \neq 0$} \nonumber
\end{eqnarray}
Below, we compute $u$ satisfying $(\ref{41})$ with each of the above cases.

\subsubsection{\textup{The cases of $a_{2,i},b_i,c=0$}} In this case, since we have $(\ref{41})$, we have
\begin{eqnarray}\label{1.1}
\sum_{i=1}^{N} a_{1,i} x_i y_1= 0 
\end{eqnarray}
The necessary and sufficient conditions to satisfy $(\ref{1.1})$ are to satisfy
\begin{eqnarray}\label{1.2}
a_{1,i} x_i y_1 + a_{1,j} x_j y_1 = 0
\end{eqnarray}
for $1 \leq i < j \leq N$. By solving (\ref{1.2}), we obtain
\begin{eqnarray}
(a_{1,i},a_{1,j})=\alpha_{1,i+j}(x_j,- x_i), \nonumber
\end{eqnarray}
where $\alpha_{1,i+j}$ is the element of $R_{N,2}$. Hence, the vector $u$ that satisfies $(\ref{1.2})$ is
\begin{eqnarray}\label{1.3}
u_1=\sum_{i=1}^{N} \sum_{j=1,j > i}^{N} \alpha_{1,i+j} ( x_i x_j y_1 e_i - x_i x_j y_1 e_j).
\end{eqnarray}

\subsubsection{\textup{The cases of $a_{1,i},b_i,c=0$}} In this case, since we have $(\ref{41})$, we have
\begin{eqnarray}\label{2.1}
\sum_{i=1}^{N} a_{2,i} x_i y_2= 0.
\end{eqnarray}
The necessary and sufficient conditions to satisfy $(\ref{2.1})$ are to satisfy
\begin{eqnarray}\label{2.2}
a_{2,i} x_i y_2 + a_{2,j} x_j y_2 = 0
\end{eqnarray}
for $1 \leq i < j \leq N$. By solving (\ref{2.2}), we obtain
\begin{eqnarray}
(a_{2,i},a_{2,j})=\alpha_{2,i+j}(x_j,- x_i), \nonumber
\end{eqnarray}
where $\alpha_{2,i+j}$ is the element of $R_{N,2}$. Hence, the vector $u$ that satisfies $(\ref{2.1})$ is
\begin{eqnarray}\label{2.3}
u_2=\sum_{i=1}^{N} \sum_{j=1,j > i}^{N} \alpha_{2,i+j} ( x_i x_j y_2 e_{N+i} - x_i x_j y_2 e_{N+j}).
\end{eqnarray}

\subsubsection{\textup{The cases of $a_{1,i},a_{2,i},c=0$}}

In this case, since we have $(\ref{41})$, we have
\begin{eqnarray}\label{3.1}
\sum_{i=2}^{N} b_i x_i s t= 0.
\end{eqnarray}
The necessary and sufficient conditions to satisfy $(\ref{3.1})$ are to satisfy
\begin{eqnarray}\label{3.2}
b_i x_i s t+ b_j x_j s t = 0
\end{eqnarray}
for $2 \leq i < j \leq N$.  By solving (\ref{3.2}), we obtain
\begin{eqnarray}
(b_i,b_j)=\alpha_{3,i+j}(x_j,- x_i), \nonumber
\end{eqnarray}
where $\alpha_{3,i+j}$ is the element of $R_{N,2}$. Hence, the vector $u$ that satisfies $(\ref{3.1})$ is
\begin{eqnarray}\label{3.3}
u_3=\sum_{i=2}^{N} \sum_{j=3,j > i}^{N} \alpha_{3,i+j} ( x_i x_j s t e_{2N+i-1} - x_i x_j s t e_{2N+j-1}).
\end{eqnarray}

\subsubsection{The cases of $b_i,c=0$}

In this case, since we have $(\ref{41})$, we have
\begin{eqnarray}\label{4.1}
\sum_{i=1}^{N} a_{1,i} x_i y_1 + \sum_{i=1}^{N} a_{2,i} x_i y_2 = 0 
\end{eqnarray}
Since we have $(\ref{1.3})$ and $(\ref{2.3})$, we only have to compute $(a_{1,i},a_{2,j})$ satisfying
\begin{eqnarray}\label{4.2}
a_{1,i} x_i y_1 + a_{2,j} x_j y_2 = 0
\end{eqnarray}
for $1 \leq i , j \leq N$. The computational result is as follows:
\begin{eqnarray}\label{4.3}
(a_{1,i},a_{2,j})=
\begin{cases}
\alpha_{4,2i}(y_2 , - y_1), & i = j, \nonumber \\
\alpha_{4,i+j}(x_j y_2 , - x_i y_1), & i \neq j, \nonumber
\end{cases}
\end{eqnarray}
where $\alpha_{4,i+j}$ is the polynomial of $R_{N,2}$ for $1 \leq i , j \leq N$.
Hence, the vector $u$ that satisfies $(\ref{5.1})$ is $u_1+u_2+u_4$, where
\begin{eqnarray}\label{4.4}
u_4
&=&
\sum_{i=1}^{N} \alpha_{4,2i} (x_i y_1 y_2 e_i - x_i y_1 y_2 e_{N+i})
+
\sum_{i=1}^{N} \sum_{j=1,i \neq j}^{N} \alpha_{4,i+j} (x_i x_j y_1 y_2 e_i  - x_i x_j y_1 y_2 e_{N+j}).
\end{eqnarray}
 We can rearrange $(\ref{4.4})$ by setting $\alpha_{4,i+j}=0$ for $i \neq j$:
\begin{eqnarray}
u_4
&=&
\sum_{i=1}^{N} \alpha_{4,2i} (x_i y_1 y_2 e_i - x_i y_1 y_2 e_{N+i})
+
\sum_{i=1}^{N} \sum_{j=1,i \neq j}^{N} \alpha_{4,i+j} (x_i x_j y_1 y_2 e_i  - x_i x_j y_1 y_2 e_{N+j})
\nonumber \\
&=&
\sum_{i=1}^{N} \Biggr ( \alpha_{4,2i} + \sum_{j=1,i \neq j}^{N} \alpha_{4,i+j} x_j \Biggr ) x_i y_1 y_2 e_i
+
\sum_{i=1}^{N} \Biggr ( - \alpha_{4,2i} x_i y_1 y_2 e_{N+i} - \sum_{j=1,i \neq j}^{N} \alpha_{4,i+j} x_i x_j y_1 y_2 e_{N+j} \Biggr )
\nonumber \\
&=&
\sum_{i=1}^{N} \alpha_{4,2i} ( x_i y_1 y_2  e_i - x_i y_1 y_2 e_{N+i}).
\nonumber 
\end{eqnarray}
Therefore, we have
\begin{eqnarray}\label{4.5}
u_4
=
\sum_{i=1}^{N} \alpha_{4,2i} ( x_i y_1 y_2  e_i - x_i y_1 y_2 e_{N+i}).
\end{eqnarray}

\subsubsection{The case of $a_{2,i},c=0$}

In this case, since we have $(\ref{41})$, we have
\begin{eqnarray}\label{5.1}
\sum_{i=1}^{N} a_{1,i} x_i y_1 + \sum_{i=2}^{N} b_i x_i s t = 0 
\end{eqnarray}
Since we have $(\ref{1.3})$ and $(\ref{3.3})$, we only have to compute $(a_{1,i},b_j)$ satisfying
\begin{eqnarray}\label{5.2}
a_{1,i} x_i y_1 + b_j x_j s t = 0
\end{eqnarray}
for $1 \leq i \leq N$ and $2 \leq j \leq N$.
The computational result is as follows:
\begin{eqnarray}\label{5.3}
(a_{1,i},b_j)=
\begin{cases}
\alpha_{5,2i}(s t , - y_1), & i = j, \nonumber \\
\alpha_{5,i+j}(x_j s t , - x_i y_1), & i \neq j, \nonumber
\end{cases}
\end{eqnarray}
where $\alpha_{5,i+j}$ is the element of $R_{N,2}$ for $1 \leq i , j \leq N$.
Hence, the vector $u$ that satisfies $(\ref{5.1})$ is $u_1+u_3+u_5$, where
\begin{eqnarray}\label{5.4}
u_5=\sum_{i=2}^{N} \alpha_{5,2i} (x_i y_1 s t e_i - x_i y_1 s t e_{2N+i-1})+
\sum_{i=1}^{N} \sum_{j=1,i \neq j}^{N} \alpha_{5,i+j} (x_i x_j y_1 s t e_i - x_i x_j y_1 s t e_{2N+j-1}).
\end{eqnarray}
We can rearrange (\ref{5.4}) by setting $\alpha_{5,i+j}=0$ for $i \neq j$:
\begin{eqnarray}
u_5&=&\sum_{i=2}^{N} \alpha_{5,2i} (x_i y_1 s t e_i - x_i y_1 s t e_{2N+i-1})+
\sum_{i=1}^{N} \sum_{j=1,i \neq j}^{N} \alpha_{5,i+j} (x_i x_j y_1 s t e_i - x_i x_j y_1 s t e_{2N+j-1})
\nonumber \\
&=&
\sum_{i=2}^{N} \alpha_{5,2i} (x_i y_1 s t e_i - x_i y_1 s t e_{2N+i-1})
+
\sum_{j=2}^{N} \alpha_{5,1+j} (x_1 x_j y_1 s t e_1 - x_1 x_j y_1 s t e_{2N+j-1})
\nonumber \\
&&
+
\sum_{i=2}^{N} \sum_{j=1,i \neq j}^{N} \alpha_{5,i+j} (x_i x_j y_1 s t e_i - x_i x_j y_1 s t e_{2N+j-1})
\nonumber \\
&=&
\sum_{i=2}^{N} \Biggr ( \alpha_{5,2i} e_i + \alpha_{5,1+j} x_1 x_j y_1 s t e_1 \Biggr ) 
-
\sum_{i=2}^{N} ( \alpha_{5,2i}  - \alpha_{5,1+j} x_i ) x_1 y_1 s t e_{2N+i-1}
\nonumber \\
&&
+
\sum_{i=2}^{N} \sum_{j=1,i \neq j}^{N} \alpha_{5,i+j} (x_i x_j y_1 s t e_i - x_i x_j y_1 s t e_{2N+j-1})
\nonumber \\
&=&
\sum_{i=2}^{N} \alpha_{5,2i} (x_i y_1 s t e_i - x_i y_1 s t e_{2N+i-1})
+
\sum_{i=2}^{N} \sum_{j=1,i \neq j}^{N} \alpha_{5,i+j} (x_i x_j y_1 s t e_i - x_i x_j y_1 s t e_{2N+j-1})
\nonumber \\
&=&
\sum_{i=2}^{N} \Biggr ( \alpha_{5,2i} + \sum_{j=1,i \neq j}^{N} \alpha_{5,i+j} x_j \Biggr ) x_i y_1 s t e_i
+
\sum_{i=2}^{N} \Biggr ( - x_i y_1 s t e_{2N+i-1} - \sum_{j=1,i \neq j}^{N} x_i x_j y_1 s t e_{2N+j-1} \Biggr )
\nonumber \\
&=&
\sum_{i=2}^{N} \alpha_{5,2i} (x_i y_1 s t e_i - x_i y_1 s t e_{2N+i-1}).
\nonumber 
\end{eqnarray}
Therefore, we have
\begin{eqnarray}\label{5.5}
u_5=\sum_{i=2}^{N} \alpha_{5,2i} (x_i y_1 s t e_i - x_i y_1 s t e_{2N+i-1}).
\end{eqnarray}

\subsubsection{\textup{The case of $a_{2,i},b_i=0,c \neq 0$}}
In this case, since we have $(\ref{41})$, we have
\begin{eqnarray}\label{6.1}
\sum_{i=1}^{N} a_{1,i} x_i y_1 + c y_2 s t = 0 .
\end{eqnarray}
The necessary and sufficient conditions to satisfy $(\ref{6.1})$ are to satisfy
\begin{eqnarray}\label{6.2}
a_{1,i} x_i y_1 + c y_2 s t = 0
\end{eqnarray}
for $1 \leq i \leq N$. By solving (\ref{6.2}), we obtain
\begin{eqnarray}
(a_{1,i},c)=\alpha_{6,i} ( y_2 s t , - x_i y_1 ), \nonumber
\end{eqnarray}
where $\alpha_{6,i}$ is the element of $R_{N,2}$ for $1 \leq i \leq N$. Hence, the vector $u$ that satisfies $(\ref{6.1})$ is
\begin{eqnarray}\label{6.3}
u_6=\sum_{i=1}^{N} \alpha_{6,i} (x_i y_1 y_2 s t e_i - x_i y_1 y_2 s t e_{3N}).
\end{eqnarray}

\subsubsection{\textup{The case of $a_{1,i},c=0$}}
In this case, since we have $(\ref{41})$, we have
\begin{eqnarray}\label{7.1}
\sum_{i=1}^{N} a_{2,i} x_i y_2 + \sum_{i=2}^{N} b_i x_i s t= 0.
\end{eqnarray}
Since we have $(\ref{2.3})$ and $(\ref{3.3})$, we only have to compute $(a_{2,i},b_i)$ satisfying
\begin{eqnarray}\label{7.2}
a_{2,i} x_i y_2 + b_j x_j s t = 0
\end{eqnarray}
for $1 \leq i , j \leq N$ and $2 \leq j \leq N$.
The computational result is as follows:
\begin{eqnarray}
(a_{2,i},b_j)=
\begin{cases}
\alpha_{7,2i}(s t , - y_2), & i = j, \nonumber \\
\alpha_{7,i+j}(x_j s t , - x_i y_2), & i \neq j. \nonumber
\end{cases}
\end{eqnarray}
where $\alpha_{7,i+j}$ is the element of $R_{N,2}$ for $1 \leq i , j \leq N$. Hence, the vector $u$ that satisfies $(\ref{7.1})$ is $u_2+u_3+u_7$, where
\begin{eqnarray}\label{7.3}
u_7=
\sum_{i=2}^{N} \alpha_{7,2i} (x_i y_2 s t e_{N+i} - x_i y_2 s t e_{2N+i-1})
+
\sum_{i=1}^{N} \sum_{j=1,i \neq j}^{N} \alpha_{7,i+j} (x_i x_j y_2 s t e_{N+i} - x_i x_j y_2 s t e_{2N+j-1}).
\end{eqnarray}
We can rearrange (\ref{7.3}) by setting $\alpha_{7,i+j}=0$ for $i \neq j$:
\begin{eqnarray}
u_7
&=&
\sum_{i=2}^{N} \alpha_{7,2i} (x_i y_2 s t e_{N+i} - x_i y_2 s t e_{2N+i-1})
+
\sum_{i=1}^{N} \sum_{j=1,i \neq j}^{N} \alpha_{7,i+j} (x_i x_j y_2 s t e_{N+i} - x_i x_j y_2 s t e_{2N+j-1})
\nonumber \\
&=&
\sum_{i=2}^{N} \alpha_{7,2i} (x_i y_2 s t e_{N+i} - x_i y_2 s t e_{2N+j-1})
+
\sum_{i=2}^{N} \alpha_{7,1+i} (x_1 x_i y_2 s t e_{N+1} - x_1 x_i y_2 s t e_{2N+i-1})
\nonumber \\
&&
+
\sum_{i=2}^{N} \sum_{j=1,i \neq j}^{N} \alpha_{7,i+j} (x_i x_j y_2 s t e_{N+i} - x_i x_j y_2 s t e_{2N+j-1})
\nonumber \\
&=&
\sum_{i=2}^{N} \alpha_{7,2i} ( x_i y_2 s t e_{N+i} + \alpha_{7,1+i} x_1 x_i y_2 s t e_{N+1} )
- 
\sum_{i=2}^{N} ( \alpha_{7,2i} + \alpha_{7,1+i} x_i ) x_1 y_2 s t e_{2N+i-1}
\nonumber \\
&&
+
\sum_{i=2}^{N} \sum_{j=1,i \neq j}^{N} \alpha_{7,i+j} (x_i x_j y_2 s t e_{N+i} - x_i x_j y_2 s t e_{2N+j-1})
\nonumber \\
&=&
\sum_{i=2}^{N} \alpha_{7,2i} ( x_i y_2 s t e_{N+i} - x_1 y_2 s t e_{2N+i-1} )
+
\sum_{i=2}^{N} \sum_{j=1,i \neq j}^{N} \alpha_{7,i+j} (x_i x_j y_2 s t e_{N+i} - x_i x_j y_2 s t e_{2N+j-1})
\nonumber \\
&=&
\sum_{i=2}^{N} 
\Biggr ( 
\biggr (
\alpha_{7,2i} 
+
\sum_{j=1,i \neq j}^{N} \alpha_{7,i+j} x_j
\biggr )
x_i y_2 s t e_{N+i}
-
\biggr (
\alpha_{7,2i}
+
\sum_{j=1,i \neq j}^{N} \alpha_{7,i+j} x_j
\biggr )
x_i y_2 s t e_{2N+j-1}
\Biggr )
\nonumber \\
&=&
\sum_{i=2}^{N} \alpha_{7,2i} (x_i y_2 s t e_{N+i} - x_i y_2 s t e_{2N+i-1}).
\nonumber
\end{eqnarray}
Therefore, we have
\begin{eqnarray}\label{7.4}
u_7
=
\sum_{i=2}^{N} \alpha_{7,2i} (x_i y_2 s t e_{N+i} - x_i y_2 s t e_{2N+i-1}).
\end{eqnarray}

\subsubsection{\textup{The case of $a_{1,i},b_i=0,c \neq 0$}}

In this case, since we have $(\ref{41})$, we have
\begin{eqnarray}\label{8.1}
\sum_{i=1}^{N} a_{2,i} x_i y_2 + c y_2 s t = 0 .
\end{eqnarray}
The necessary and sufficient conditions to satisfy $(\ref{8.1})$ are to satisfy
\begin{eqnarray}\label{8.2}
a_{2,i} x_i y_2 + c y_2 s t = 0
\end{eqnarray}
for $1 \leq i \leq N$. By solving (\ref{8.2}), we obtain
\begin{eqnarray}
(a_{2,i},c)=\alpha_{8,i} ( s t , - x_i). \nonumber
\end{eqnarray}
where $\alpha_{8,i}$ is the element of $R_{N,2}$ for $1 \leq i \leq N$. Hence, the vector $u$ that satisfies $(\ref{8.1})$ is
\begin{eqnarray}\label{8.3}
u_8=\alpha_{8,1} (x_1 y_2 s t e_{N+i} - x_1 y_2 s t e_{3N}) + \sum_{i=2}^{N} \alpha_{8,i} (x_i y_2 s t e_{N+i} - x_i y_2 s t e_{3N}).
\end{eqnarray}

\subsubsection{\textup{The case of $a_{1,i},a_{2,i}=0,c \neq 0$}}
In this case, since we have $(\ref{41})$, we have
\begin{eqnarray}\label{9.1}
\sum_{i=2}^{N} b_i x_i s t  + c y_2 s t = 0 
\end{eqnarray}
The necessary and sufficient conditions to meet $(\ref{9.1})$ are
\begin{eqnarray}\label{9.2}
b_i x_i s t  + c y_2 s t = 0
\end{eqnarray}
for $2 \leq i \leq N$. By solving (\ref{9.2}), we obtain
\begin{eqnarray}
(b_i,c)=\alpha_{9,i} (y_2,x_i)
\end{eqnarray}
where $\alpha_{9,i}$ is the element of $R_{N,2}$. Hence, the vector $u$ that satisfies $(\ref{9.1})$ is
\begin{eqnarray}\label{9.3}
u_9=\sum_{i=2}^{N} \alpha_{9,i} ( x_i y_2 s t e_{2N+i-1} - x_i y_2 s t e_{3N} ).
\end{eqnarray}

\subsection{The computational of component of $u$} \label{sub2}
By (\ref{1.3}), (\ref{2.3}), (\ref{3.3}), (\ref{4.5}), (\ref{5.5}), (\ref{6.3}), (\ref{7.4}), (\ref{8.3}) and (\ref{9.3}), we have $u=\sum_{i=1}^9 u_i$. Hence, the positive component $u^{+}$ of $u$ is as follows:
\begin{eqnarray}
u^{+}
&=&
\sum_{i=1}^{N} \sum_{j=1,j > i}^{N} \alpha_{1,i+j} x_i x_j y_1 e_i 
+
\sum_{i=1}^{N} \sum_{j=1,j > i}^{N} \alpha_{2,i+j} x_i x_j y_2 e_{N+i}
+
\sum_{i=2}^{N} \sum_{j=1,j > i}^{N} \alpha_{3,i+j} x_i x_j s t e_{2N+i-1}
\nonumber \\
&&
+
\sum_{i=1}^{N} \alpha_{4,2i} x_i y_1 y_2 e_i
+
\sum_{i=1}^{N} \alpha_{5,2i} x_i y_1 s t e_i
+
\sum_{i=1}^{N} \alpha_{6,i} x_i y_1 y_2 s t e_i
\nonumber \\
&&
+
\sum_{i=2}^{N} \alpha_{7,2i} x_i y_2 s t e_{N+i}
+
\alpha_{8,1} x_1 y_2 s t e_{N+1}
+
\sum_{i=2}^{N} \alpha_{8,i} x_i y_2 s t e_{N+i}
+
\sum_{i=2}^{N} \alpha_{9,i} x_i y_2 s t e_{2N+i-1}
\nonumber
\end{eqnarray}
Moreover, the negative component $u^{-}$ of $u$ is as follows:
\begin{eqnarray}
u^{-}
&=&
\sum_{i=1}^{N} \sum_{j=1,j > i}^{N} \alpha_{1,i+j} x_i x_j y_1 e_j
+
\sum_{i=1}^{N} \sum_{j=1,j > i}^{N} \alpha_{2,i+j} x_i x_j y_2 e_{N+j}
+
\sum_{i=1}^{N} \sum_{j=1,j > i}^{N} \alpha_{3,i+j} x_i x_j s t e_{2N+j-1}
\nonumber \\
&&
\sum_{i=1}^{N} \alpha_{4,2i}  x_i y_1 y_2 e_{N+i}
+
\sum_{i=2}^{N} \alpha_{5,2i} x_i y_1 s t e_{2N+i-1}
+
\sum_{i=1}^{N} \alpha_{6,i} x_i y_1 y_2 s t e_{3N}
\nonumber \\
&&
+
\sum_{i=2}^{N} \alpha_{7,2i} x_i y_2 s t e_{2N+i-1}
+
\alpha_{8,1} x_1 y_2 s t e_{3N}
+
\sum_{i=2}^{N} \alpha_{8,i} x_i y_2 s t e_{3N}
+
\sum_{i=1}^{N} \alpha_{9,i} x_i y_2 s t e_{3N}
\nonumber
\end{eqnarray}

Then, we compute the component of $u=u^{+}-u^{-}$ in detail.

\subsubsection{\textup{The computational of the $i$-th component of u for $1 \leq i \leq N$}}

For $1 \leq i \leq N$, the $i$-th positive component of u is is computed the following:
\begin{eqnarray}
&&
\sum_{i=1}^{N} \sum_{j=1,j > i}^{N} \alpha_{1,i+j} x_i x_j y_1 e_i
+
\sum_{i=1}^{N} \alpha_{4,2i} x_i y_1 y_2 e_i
+
\sum_{i=2}^{N} \alpha_{5,2i} x_i y_1 s t e_i
+
\sum_{i=1}^{N} \alpha_{6,i} x_i y_1 y_2 s t e_i
\nonumber \\
&=&
\sum_{j=2}^{N} \alpha_{1,1+j} x_1 x_j y_1 e_1
+
\alpha_{4,2} x_1 y_1 y_2 e_1
+
\sum_{i=2}^{N}
\Biggr (
\sum_{j=1,j > i}^{N} \alpha_{1,i+j} x_j
+
\alpha_{4,2i} y_2
+
\alpha_{5,2i} s t
+
\alpha_{6,i} y_2 s t 
\Biggr ) x_i y_1 e_i
\nonumber \\
&=&
\sum_{j=2}^{N} \alpha_{1,1+j} x_1 x_j y_1 e_1
+
\alpha_{4,2} x_1 y_1 y_2 e_1
+
\sum_{i=1}^{N}
\Biggr (
\sum_{j=1,j > i}^{N} \alpha_{1,i+j} x_j
+
( 
\alpha_{4,2i}
+
\alpha_{6,i} s t
)
y_2 
+
\alpha_{5,2i} s t
\Biggr ) x_i y_1 e_i
\nonumber \\
&=&
\sum_{j=2}^{N} \alpha_{1,1+j} x_1 x_j y_1 e_1
+
\alpha_{4,2} x_1 y_1 y_2 e_1
+
\sum_{i=2}^{N}
\Biggr (
\sum_{j=1,j > i}^{N} \alpha_{1,i+j} x_j
+ 
\alpha_{4,2i} y_2 
+
\alpha_{5,2i} s t
\Biggr ) x_i y_1 e_i.
\nonumber \\
&=&
\sum_{j=2}^{N} \alpha_{1,1+j} x_1 x_j y_1 e_1
+
\alpha_{4,2} x_1 y_1 y_2 e_1
+
\sum_{i=2}^{N} \sum_{j=1,j > i}^{N} \alpha_{1,i+j} x_i x_j y_1 e_i
+ 
\sum_{i=2}^{N} \alpha_{4,2i} x_i y_1 y_2 e_i
+
\sum_{i=2}^{N} \alpha_{5,2i} x_i y_1 s t e_i
\nonumber
\end{eqnarray}
On the other hand, for $1 \leq i \leq N$, the $i$-th negative component of u is 
\begin{eqnarray}
\sum_{i=1}^{N} \sum_{j=1,j > i}^{N} \alpha_{1,i+j} x_i x_j y_1 e_j=\sum_{j=2}^{N} \alpha_{1,1+j} x_1 x_j y_1 e_j+ \sum_{i=2}^{N} \sum_{j=1,j > i}^{N} \alpha_{1,i+j} x_i x_j y_1 e_j. \nonumber
\end{eqnarray}
Hence, for $1 \leq i \leq N$, the $i$-th component of u is 
\begin{eqnarray}\label{component1}
&&
\sum_{j=2}^{N} \alpha_{1,1+j} x_1 x_j y_1 e_1
+
\alpha_{4,2} x_1 y_1 y_2 e_1
+
\sum_{i=2}^{N} \sum_{j=1,j > i}^{N} \alpha_{1,i+j} x_i x_j y_1 e_i
+ 
\sum_{i=2}^{N} \alpha_{4,2i} x_i y_1 y_2 e_i
+
\sum_{i=2}^{N} \alpha_{5,2i} x_i y_1 s t e_i
\nonumber \\
&&
-
\sum_{j=2}^{N} \alpha_{1,1+j} x_1 x_j y_1 e_j
-
\sum_{i=2}^{N} \sum_{j=1,j > i}^{N} \alpha_{1,i+j} x_i x_j y_1 e_j.
\end{eqnarray}

\subsubsection{\textup{The computational of the $i$-th positive component of u for $N+1 \leq i \leq 2N$}}

For $1 \leq i \leq N$, the $N+i$-th positive component of u is as follows:
\begin{eqnarray}
\sum_{i=1}^{N} \sum_{j=1,j > i}^{N} \alpha_{2,i+j} x_i x_j y_2 e_{N+i}
+
+
\sum_{i=2}^{N} \alpha_{7,2i} x_i y_2 s t e_{N+i}
+
\alpha_{8,1} x_1 y_2 s t e_{N+1}
+
\sum_{i=2}^{N} \alpha_{8,i} x_i y_2 s t e_{N+i}
\nonumber
\end{eqnarray}
On the other hand, for $1 \leq i \leq N$, the $i$-th negative component of u is is as follows:
\begin{eqnarray}
\sum_{i=1}^{N} \sum_{j=1,j > i}^{N} \alpha_{2,i+j} x_i x_j y_2 e_{N+j}
+
\sum_{i=1}^{N} \alpha_{4,2i}  x_i y_1 y_2 e_{N+i}.
\nonumber
\end{eqnarray}
Hence, for $N+1 \leq i \leq 2N$, the $i$-th component of u is 
\begin{eqnarray}\label{component2}
&&
\sum_{i=1}^{N} \sum_{j=1,j > i}^{N} \alpha_{2,i+j} x_i x_j y_2 e_{N+i}
+
\sum_{i=2}^{N} \alpha_{7,2i} x_i y_2 s t e_{N+i}
+
\alpha_{8,1} x_1 y_2 s t e_{N+1}
+
\sum_{i=2}^{N} \alpha_{8,i} x_i y_2 s t e_{N+i}
\nonumber \\
&&
\quad
-
\sum_{i=1}^{N} \sum_{j=1,j > i}^{N} \alpha_{2,i+j} x_i x_j y_2 e_{N+j}
-
\sum_{i=1}^{N} \alpha_{4,2i}  x_i y_1 y_2 e_{N+i}.
\end{eqnarray}

\subsubsection{\textup{The computational of the $i$-th positive component of u for $2N+1 \leq i \leq 3N$}}
First, we compute the $i$-th positive component of u for $2N+1 \leq i \leq 3N-1$. For $2 \leq i \leq N-1$, the $2N+i-1$-th positive component of u is as follows following:
\begin{eqnarray}
\sum_{i=2}^{N} \sum_{j=1,j > i}^{N} \alpha_{3,i+j} x_i x_j s t e_{2N+i-1}
+
\sum_{i=2}^{N} \alpha_{9,i} x_i y_2 s t e_{2N+i-1}
\nonumber
\end{eqnarray}
On the other hand, for $1 \leq i \leq N-1$, the $2N+i-1$-th negative component of u is as follows:
\begin{eqnarray}
\sum_{i=2}^{N} \sum_{j=1,j > i}^{N} \alpha_{3,i+j} x_i x_j s t e_{2N+j-1}
+
\sum_{i=2}^{N} \alpha_{5,2i} x_i y_1 s t e_{2N+i-1}
+
\sum_{i=2}^{N} \alpha_{7,2i} x_i y_2 s t e_{2N+i-1}
\nonumber
\end{eqnarray}
Hence, for $2N+1 \leq i \leq 3N$, the computational result of the $i$-th component of u is as follows:
\begin{eqnarray}\label{component3}
&&
\sum_{i=2}^{N} \sum_{j=1,j > i}^{N} \alpha_{3,i+j} x_i x_j s t e_{2N+i-1}
+
\sum_{i=2}^{N} \alpha_{9,i} x_i y_2 s t e_{2N+i-1}
\nonumber \\
&&
\quad
-
\sum_{i=2}^{N} \sum_{j=1,j > i}^{N} \alpha_{3,i+j} x_i x_j s t e_{2N+j-1}
-
\sum_{i=2}^{N} \alpha_{5,2i} x_i y_1 s t e_{2N+i-1}
-
\sum_{i=2}^{N} \alpha_{7,2i} x_i y_2 s t e_{2N+i-1}
\end{eqnarray}

Second, the $3N$-th component of u is computed by using $\alpha_{6,i}=0$ as follows:
\begin{eqnarray}
&&
-
\sum_{i=1}^{N} \alpha_{6,i} x_i y_1 y_2 s t e_{3N}
-
\alpha_{8,1} x_1 y_2 s t e_{3N}
-
\sum_{i=2}^{N} \alpha_{8,i} x_i y_2 s t e_{3N}
-
\sum_{i=2}^{N} \alpha_{9,i} x_i y_2 s t e_{3N}
\nonumber \\
&=&
\label{component4}
-
\alpha_{8,1} x_1 y_2 s t e_{3N}
-
\sum_{i=2}^{N} \alpha_{8,i} x_i y_2 s t e_{3N}
-
\sum_{i=2}^{N} \alpha_{9,i} x_i y_2 s t e_{3N}
\end{eqnarray}

\subsection{Matrix representation of $u$}

By (\ref{component1}), (\ref{component2}), (\ref{component3}) and (\ref{component4}), we have
\begin{eqnarray}\label{component5}
u
&=&
\sum_{j=2}^{N} \alpha_{1,1+j} x_1 x_j y_1 e_1
+
\alpha_{4,2} x_1 y_1 y_2 e_1
+
\sum_{i=2}^{N} \sum_{j=1,j > i}^{N} \alpha_{1,i+j} x_i x_j y_1 e_i
+ 
\sum_{i=2}^{N} \alpha_{4,2i} x_i y_1 y_2 e_i
+
\sum_{i=2}^{N} \alpha_{5,2i} x_i y_1 s t e_i
\nonumber \\
&&
-
\sum_{j=2}^{N} \alpha_{1,1+j} x_1 x_j y_1 e_j
-
\sum_{i=2}^{N} \sum_{j=1,j > i}^{N} \alpha_{1,i+j} x_i x_j y_1 e_j
\nonumber \\
&&
+
\sum_{i=1}^{N} \sum_{j=1,j > i}^{N} \alpha_{2,i+j} x_i x_j y_2 e_{N+i}
+
\sum_{i=2}^{N} \alpha_{7,2i} x_i y_2 s t e_{N+i}
+
\alpha_{8,1} x_1 y_2 s t e_{N+1}
+
\sum_{i=2}^{N} \alpha_{8,i} x_i y_2 s t e_{N+i}
\nonumber \\
&&
\quad
-
\sum_{i=1}^{N} \sum_{j=1,j > i}^{N} \alpha_{2,i+j} x_i x_j y_2 e_{N+j}
-
\sum_{i=1}^{N} \alpha_{4,2i}  x_i y_1 y_2 e_{N+i}
\nonumber \\
&&
+
\sum_{i=2}^{N} \sum_{j=1,j > i}^{N} \alpha_{3,i+j} x_i x_j s t e_{2N+i-1}
+
\sum_{i=2}^{N} \alpha_{9,i} x_i y_2 s t e_{2N+i-1}
\nonumber \\
&&
\quad
-
\sum_{i=2}^{N} \sum_{j=1,j > i}^{N} \alpha_{3,i+j} x_i x_j s t e_{2N+j-1}
-
\sum_{i=2}^{N} \alpha_{5,2i} x_i y_1 s t e_{2N+i-1}
-
\sum_{i=2}^{N} \alpha_{7,2i} x_i y_2 s t e_{2N+i-1}
\nonumber \\
&&
-
\alpha_{8,1} x_1 y_2 s t e_{3N}
-
\sum_{i=2}^{N} \alpha_{8,i} x_i y_2 s t e_{3N}
-
\sum_{i=2}^{N} \alpha_{9,i} x_i y_2 s t e_{3N}.
\end{eqnarray}
By using new vector $d_{i,j}=e_i - e_j$ for $1 \leq i , j \leq 3N$, we can rearrange (\ref{component5}) as follows:
\begin{eqnarray}\label{component6}
u
&=&
\sum_{j=2}^{N} d_{1,j} \alpha_{1,1+j} x_1 x_j y_1
+
\sum_{i=2}^{N} \sum_{j=1,j > i}^{N} d_{i,j} \alpha_{1,i+j} x_i x_j y_1
+
\sum_{i=1}^{N} \sum_{j=1,j > i}^{N} d_{N+i,N+j} \alpha_{2,i+j} x_i x_j y_2
\nonumber \\
&&
+
\sum_{i=1}^{N} d_{i,N+i} \alpha_{4,2i}  x_i y_1 y_2
+
\sum_{i=2}^{N} \sum_{j=1,j > i}^{N} d_{2N+i-1,2N+j-1} \alpha_{3,i+j} x_i x_j s t 
+
\sum_{i=2}^{N} d_{i,2N+i-1} \alpha_{5,2i} x_i y_1 s t
\nonumber \\
&&
+
d_{N+1,3N} \alpha_{8,1} x_1 y_2 s t
+
\sum_{i=2}^{N} d_{N+i,3N} \alpha_{8,i} x_i y_2 s t
+
\sum_{i=2}^{N} d_{2N+i-1,3N} \alpha_{9,i} x_i y_2 s t .
\end{eqnarray}
We introduce the symbols necessary to express u by using a matrix. Let denote $g_1(N)$ and $g_2(N)$ by
\begin{eqnarray}
g_1(N)
&=&
v_{1,1} \boxplus v_{1,2} \boxplus \cdots \boxplus v_{1,N-1} \boxplus v_{2,1} \boxplus v_{2,2} \boxplus \cdots \boxplus v_{2,N-1} \boxplus v_{3,N},
\nonumber \\
g_2(N)
&=&
v_{4,2} \boxplus v_{4,3} \boxplus \cdots \boxplus v_{4,N-1} \boxplus v_{5,N} \boxplus v_{6,N} \boxplus v_{7,N},
\nonumber
\end{eqnarray}
where
\begin{eqnarray}
v_{1,i}
&=&
\begin{pmatrix}
x_i x_{i+1} y_1 & x_i x_{i+2} y_1 & \cdots & x_i x_N y_1 \\
\end{pmatrix}
, \quad \text{for $1 \leq i \leq N-1$},
\nonumber \\
v_{2,i}
&=&
\begin{pmatrix}
x_i x_{i+1} y_2 & x_i x_{i+2} y_2 & \cdots & x_i x_N y_1 \\
\end{pmatrix}
, \quad \text{for $1 \leq i \leq N-1$},
\nonumber \\
v_{3,N}
&=&
\begin{pmatrix}
x_1 y_1 y_2 & x_2 y_1 y_2 & \cdots & x_N y_1 y_2 \\
\end{pmatrix}
,
\nonumber \\
v_{4,i}
&=&
\begin{pmatrix}
x_i x_{i+1} s t & x_i x_{i+2} s t & \cdots & x_i x_N s t \\
\end{pmatrix}
, \quad \text{for $2 \leq i \leq N-1$},
\nonumber \\
v_{5,N}
&=&
\begin{pmatrix}
x_2 y_1 s t & x_3 y_1 s t & \cdots & x_N y_1 s t \\
\end{pmatrix}
,
\nonumber \\
v_{6,N}
&=&
\begin{pmatrix}
x_1 y_2 s t & x_2 y_2 s t & x_3 y_2 s t & \cdots & x_N y_2 s t \\
\end{pmatrix}
,
\nonumber \\
v_{7,N}
&=&
\begin{pmatrix}
x_2 y_2 s t & x_3 y_2 s t & \cdots & x_N y_2 s t \\
\end{pmatrix}
\nonumber
.
\nonumber
\end{eqnarray}
Let denote the matrix $f_{N,2}$ by
\begin{eqnarray}
&&
\Biggr ( \bigotimes_{i=2}^N d_{1,i} \Biggr )
\otimes
\Biggr ( \bigotimes_{i=2}^{N} \bigotimes_{j=3,j > i}^{N} d_{i,j} \Biggr )
\otimes
\Biggr ( \bigotimes_{i=1}^{N} \bigotimes_{j=2,j > i}^{N} d_{N+i,N+j} \Biggr )
\otimes
\Biggr ( \bigotimes_{i=1}^{N} d_{i,N+i} \Biggr )
\nonumber \\
&&
\otimes
\Biggr ( \bigotimes_{i=2}^{N} \bigotimes_{j=1,j > i}^{N} d_{2N+i-1,2N+j-1} \Biggr )
\otimes
\Biggr ( \bigotimes_{i=2}^{N} d_{i,2N+i-1} \Biggr )
\otimes
d_{N+1,3N}
\otimes
\Biggr ( \bigotimes_{i=2}^{N} d_{N+i,3N} \Biggr )
\otimes
\Biggr ( \bigotimes_{i=2}^{N} d_{2N+i-1,3N} \Biggr ).
\nonumber
\end{eqnarray}
By using $f_{N,2}$ and $g_{N,2}:=g_1(N) \oplus g_2(N)$, it follows from (\ref{component6}) that we have $u=f_{N,2}{g_{N,2}}^{T}$. Hence, we have a free resolution
\begin{eqnarray}\label{resolution3}
0 \leftarrow R_{N,2} \xleftarrow{
\overbrace{
\begin{bmatrix}
1 & \cdots & 1
\end{bmatrix}
}^{3N \text{pieces}}
} F^{(N,2)}_{1} \xleftarrow{f_{N,2}} F'_{2} \leftarrow \cdots
\end{eqnarray}
where
\begin{eqnarray}
F'_{2}
&=&
\Biggr ( \bigoplus_{i=1}^{N} \bigoplus_{j=2,j>i}^{N} R_{N,2}(x_i x_j y_1) \Biggr )
\oplus
\Biggr ( \bigoplus_{i=1}^{N} \bigoplus_{j=2,j>i}^{N} R_{N,2}(x_i x_j y_2) \Biggr )
\oplus
\Biggr ( \bigoplus_{i=1}^{N} R_{N,2}(x_i y_1 y_2) \Biggr )
\nonumber \\
&&
\oplus
\Biggr ( \bigoplus_{i=2}^{N} \bigoplus_{j=2,j>i}^{N} R_{N,2}(x_i x_j s t) \Biggr )
\oplus
\Biggr ( \bigoplus_{i=2}^{N} R_{N,2}(x_i y_1 s t) \Biggr )
\oplus
\Biggr ( \bigoplus_{i=1}^{N} R_{N,2}(x_i y_2 s t) \Biggr )
\oplus
\Biggr ( \bigoplus_{i=2}^{N} R_{N,2}(x_i y_2 s t) \Biggr ).
\nonumber
\end{eqnarray}

\subsection{Construction of $F^{(N,2)}_{2}$}
Then, we show that $F^{(N,2)}_{2}=F'_{2}$. Suppose that $F^{(N,2)}_{2} \neq F'_{2}$. Since we have $F^{(N,2)}_{2} \subsetneq F'_{2}$, for instance we can suppose that
\begin{eqnarray}
F^{(N,2)}_{2}
&=&
\Biggr ( \bigoplus_{i=1}^{N} \bigoplus_{j=2,j>i}^{N} R_{N,2}(x_i x_j y_1) \Biggr )
\oplus
\Biggr ( \bigoplus_{i=1}^{N} \bigoplus_{j=2,j>i}^{N} R_{N,2}(x_i x_j y_2) \Biggr )
\oplus
\Biggr ( \bigoplus_{i=1}^{N} R_{N,2}(x_i y_1 y_2) \Biggr )
\nonumber \\
&&
\oplus
\Biggr ( \bigoplus_{i=2}^{N} \bigoplus_{j=2,j>i}^{N} R_{N,2}(x_i x_j s t) \Biggr )
\oplus
\Biggr ( \bigoplus_{i=2}^{N} R_{N,2}(x_i y_1 s t) \Biggr )
\oplus
\Biggr ( \bigoplus_{i=1}^{N} R_{N,2}(x_i y_2 s t) \Biggr ).
\nonumber
\end{eqnarray}
However, such assumption means that we have $F^{(N,2)}_{2}=0$. In fact, if we have $F^{(N,2)}_{2} \neq 0$, then we have $\im(\phi^{(N,2)}_{2}) \subsetneq \text{span}_{R_{N,2}} \{ \omega \mid \text{$\omega$ is each term of $u$.} \}=\Ker(\phi^{(N,2)}_{1})$. It contradicts that we have $\im(\phi^{(N,2)}_2)=\Ker(\phi^{(N,2)}_1)$. Since we have $F^{(N,2)}_{2}=0$, it contradicts existence of the free resolution (\ref{resolution3}). Hence, we showed that $F^{(N,2)}_{2}=F'_{2}$.

\subsection{Computational of total Betti number of $F^{(N,2)}_{2}$}
By $F^{(N,2)}_{2}=F'_{2}$, we have
\begin{eqnarray}
F^{(N,2)}_{2}
&=&
\Biggr ( \bigoplus_{i=1}^{N} \bigoplus_{j=2,j>i}^{N} R_{N,2}(x_i x_j y_1) \Biggr )
\oplus
\Biggr ( \bigoplus_{i=1}^{N} \bigoplus_{j=2,j>i}^{N} R_{N,2}(x_i x_j y_2) \Biggr )
\oplus
\Biggr ( \bigoplus_{i=1}^{N} R_{N,2}(x_i y_1 y_2) \Biggr )
\nonumber \\
&&
\oplus
\Biggr ( \bigoplus_{i=2}^{N} \bigoplus_{j=2,j>i}^{N} R_{N,2}(x_i x_j s t) \Biggr )
\oplus
\Biggr ( \bigoplus_{i=2}^{N} R_{N,2}(x_i y_1 s t) \Biggr )
\oplus
\Biggr ( \bigoplus_{i=1}^{N} R_{N,2}(x_i y_2 s t) \Biggr )
\oplus
\Biggr ( \bigoplus_{i=2}^{N} R_{N,2}(x_i y_2 s t) \Biggr ).
\nonumber
\end{eqnarray}
The above free module $F^{(N,2)}_{2}$ can be written as follows:
\begin{eqnarray}
&&
F^{(N,2)}_{2}
\nonumber \\
&=&
R_{N,2}(-3)^{{}_N C_2}
\oplus
R_{N,2}(-3)^{{}_N C_2}
\oplus
R_{N,2}(-3)^{N}
\oplus
R_{N,2}(-4)^{{}_{N-1} C_2}
\oplus
R_{N,2}(-4)^{N-1}
\oplus
R_{N,2}(-4)^N
\oplus
R_{N,2}(-4)^{N-1}
\nonumber \\
&=&
R_{N,2}(-3)^{2 {}_N C_2 + N} \oplus R_{N,2}(-4)^{{}_{N-1} C_2 + 3N-2}.
\nonumber
\end{eqnarray}
Hence, it follows from the above result that we have
\begin{eqnarray}
\sum_{i=0}^{\infty} B_{2,i}(N,2)
&=&
B_{2,3}(N,2)
+
B_{2,4}(N,2)
\nonumber \\
&=&
( 2 {}_N C_2 + N )+( {}_{N-1} C_2 + 3N-2 )
\nonumber \\
&=&
2 {}_N C_2 + {}_{N-1} C_2 + 4N-2
\nonumber \\
&=&
\frac{2N(N-1)}{2}+\frac{(N-1)(N-2)}{2}+\frac{8N-4}{2}
\nonumber \\
&=&
\frac{(2N^2-2N)+(N^2-3N+2)+(8N-4)}{2}
\nonumber \\
&=&\frac{3 N^2 + 3 N - 2}{2}
\nonumber \\
&=&N(N+1)-1.
\nonumber
\end{eqnarray}
Therefore, we showed that equation (\ref{22}) holds.

\section{Example}

In this section, we introduce computational result by using Theorem \ref{theo1} and \ref{theo2}. In particular, we mainly introduce computational result for $2 \leq n_1 \leq 20$ and $n_2=1,2$.\\

At first, for $n_2=1,2$, Computational result of total of each of the 1-th Betti number 
\begin{eqnarray}
\sum_{i=0}^{\infty} B_{1,i}(n_1,n_2)=B_{1,2}(n_1,n_2)+B_{1,3}(n_1,n_2) \nonumber
\end{eqnarray}
is as Table \ref{table1}:
\begin{table}[h]
 \caption{Total of the 1-th Betti number}
 \label{table1}
 \centering
  \begin{tabular}{|c|c|c||c|c|c|}
   \hline
   Value of $n_1$ & Value of  $B(n_1,1)$ & Value of  $B(n_1,2)$ & Value of $n_1$ & Value of  $B(n_1,1)$ & Value of  $B(n_1,2)$ \\
   \hline
   2 & 3 & 6 & 11 & 21 & 33 \\
   3 & 5 & 9 & 12 & 23 & 36 \\
   4 & 7 & 12 & 14 & 27 & 42 \\
   5 & 9 & 15 & 15 & 29 & 45 \\
   6 & 11 & 18 & 16 & 31 & 48 \\
   7 & 13 & 21 & 17 & 33 & 51 \\
   8 & 15 & 24 & 18 & 35 & 54 \\
   9 & 17 & 27 & 19 & 37 & 57 \\
   10 & 19 & 30 & 20 & 39 & 60 \\
   \hline
  \end{tabular}
\end{table}

Since we have 
\begin{eqnarray}
&&\initial_{<}(I_{L_2(n_1,1)})=\{ x_i y_1 \mid 1 \leq i \leq n_1 \} \cup \{ x_i s t \mid 2 \leq i \leq n_1 \}, \nonumber \\ 
&&\initial_{<}(I_{L_2(n_1,2)})=\{ x_i y_1 \mid 1 \leq i \leq n_1 \} \cup \{ x_i y_2 \mid 1 \leq i \leq n_1 \} \cup \{ x_i s t \mid 2 \leq i \leq n_1 \} \cup \{ y_2 s t \}\nonumber
\end{eqnarray}
by \cite{2} or \cite{4}, note that Table \ref{table1} is a trivial result.\\

\begin{figure}[h]
\centering
\includegraphics[keepaspectratio, scale=0.60]{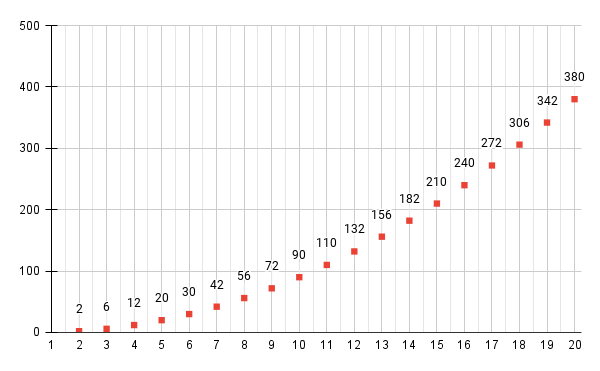}
\caption{Total of the 2-th Betti number  for $n_2=1$}
\label{chart1}
\end{figure}

Next, for total of the 2-th Betti number 
\begin{eqnarray}
\sum_{i=0}^{\infty} B_{2,i}(n_1,n_2)=B_{2,3}(n_1,n_2)+B_{2,4}(n_1,n_2), \nonumber
\end{eqnarray}
computational result in the case of $n_2=1$ is as Figure \ref{chart1}. This result was computed by using Theorem \ref{theo1}. The horizontal axis in Figure \ref{chart1} display value of $n_1$. The vertical axis in Figure \ref{chart1} display value of total of each of the 2-th Betti number $\sum_{i=0}^{\infty} B_{2,i}(n_1,1)$. The value of above the red squares in Figure \ref{chart1} display the actual calculated values.

\begin{figure}[h]
\centering
\includegraphics[keepaspectratio, scale=0.60]{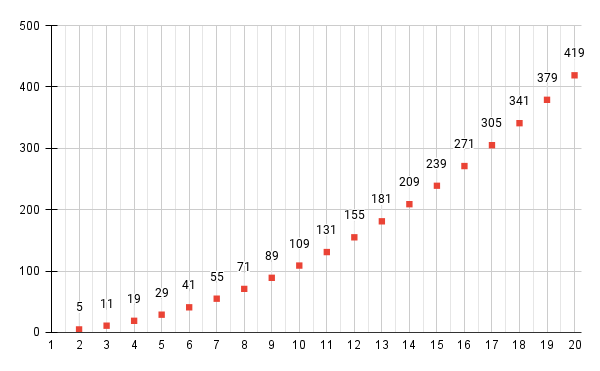}
\caption{Total of the 2-th Betti number for $n_2=2$}
\label{chart2}
\end{figure}

On the other hand, by using Theorem \ref{theo2},  computational result in the case of $n_2=2$ is as Figure \ref{chart2}. The horizontal axis in Figure \ref{chart1} display value of $n_1$. The vertical axis in Figure \ref{chart2} display value of total of each of the 2-th Betti number $\sum_{i=0}^{\infty} B_{2,i}(n_1,2)$. The value of above the red squares in Figure \ref{chart2} display the actual computed values.

\section*{Comment}

In fact, \cite{2} and \cite{4} are almost the same. The official one is \cite{2}. However, it may be better to see \cite{4} in terms of visibility.

\large

\end{document}